\documentclass[12pt,a4paper]{amsart}
\usepackage{amsmath, amsthm, amsfonts, amssymb, graphicx, hyperref, bm,verbatim,siunitx}

\theoremstyle{plain}
\newtheorem{thrm}{Theorem}
\newtheorem{lmm}[thrm]{Lemma}

\newtheorem{crllry}[thrm]{Corollary}

\newtheorem{cnjctr}{Conjecture}

\newtheorem*{dfntn}{Definition}

\numberwithin{sblmm}{thrm} 
\numberwithin{equation}{section}

\DeclareMathOperator*{\lcm}{lcm}
\DeclareMathOperator*{\Li}{Li}
\renewcommand{\phi}{\varphi}

\newcommand{\Mod}[1]{\ (\mathrm{mod}\ #1)}

\begin{document}
\title{On the Twin Prime Conjecture}
\author{James Maynard}
\email{james.alexander.maynard@gmail.com}
\begin{abstract}
We discuss various recent advances on weak forms of the Twin Prime Conjecture.
\end{abstract}
\maketitle
\section{Introduction}
One of the most famous problems in mathematics is the Twin Prime Conjecture.
\begin{cnjctr}[Twin Prime Conjecture]\label{cnjctr:Twin}
There are infinitely many pairs of primes which differ by $2$.
\end{cnjctr}
There is some debate as to how old the Twin Prime Conjecture is - it was certainly considered by de Polignac over 150 years ago \cite{Polignac}, but there has been speculation that it could go back much further, potentially as far back as Euclid and the ancient Greeks over 2000 years ago.

Unfortunately the Twin Prime Conjecture remains unsolved, but it is strongly believed to be true. Hardy and Littlewood \cite{HardyLittlewood} made a more precise conjecture about the density of twin primes. In its most optimistic form, this is the following.
\begin{cnjctr}[Quantitative Twin Prime Conjecture]\label{cnjctr:TwinPrime}
For any fixed $\epsilon>0$, the number of twin primes up to $x$, $\pi_2(x)$, is given by
\[\pi_2(x)=\#\{p\le x:\,p+2\text{ prime}\}=\mathfrak{S}\int_2^x\frac{dt}{\log^2{t}}+O(x^{1/2+\epsilon}),\]
where
\[
\mathfrak{S}=2\prod_{p\ge 3}\Bigl(1-\frac{1}{(p-1)^2}\Bigr).
\]
\end{cnjctr}
(Here, and throughout the paper, $\epsilon$ will denote a small positive quantity, and various implied constants in $O(\cdot)$ notation will be allowed to depend on $\epsilon$ without us explicitly saying so. In particular, the implied constant in the term $O(x^{1/2+\epsilon})$ may depend on $\epsilon$.)
\subsection{Numerical evidence}
The quantitative variant of the Twin Prime Conjecture given by Conjecture \ref{cnjctr:TwinPrime} allows one to test it numerically, and the evidence is as compelling as one could possibly hope for. A huge amount of computation has been done. If we let $\Li_2(x)$ denote the integral in Conjecture \ref{cnjctr:TwinPrime}, then we have the following.
\begin{center}
\begin{tabular}{|c|c|c|c|}
\hline
$x$ & $\pi_2(x)$ & $\mathfrak{S}\Li_2(x)$ & $\pi_2(x)-\mathfrak{S}\Li_2(x)$\\
\hline
10 & 2 & 4.8\dots & -2.8\dots\\
$10^2$ & 8 & 13.5\dots & -5.5\dots\\
$10^3$ & 35 & 45.7\dots & -10.7\dots\\
$10^4$ & 205 & 214.2\dots & -9.2\dots\\
$10^5$ & 1224 & 1248.7\dots & -24.7\dots\\
$10^6$ & 8169 & 8248.0\dots & -79.0\dots\\
$10^7$ & 58980 & 58753.8\dots & 226.1\dots\\
$10^8$ & 440312 & 440367.7\dots & -55.7\dots\\
$10^9$ & 3424506 & 3425308.1\dots & -802.1\dots\\
$10^{10}$ & 27412679 & 27411416.5\dots & 1262.4\dots\\
$10^{11}$ & 224376048 & 224368864.7\dots & 7183.2\dots\\
$10^{12}$ & 1870585220 & 1870559867.8\dots & 25352.1\dots\\
$10^{13}$ & 15834664872 & 15834598305.0\dots & 66566.9\dots\\
$10^{14}$ & 135780321665 & 135780264894.4\dots & 56770.5\dots\\
$10^{15}$ & 1177209242304 & 1177208491860.6\dots & 750443.3\dots\\
\hline
\end{tabular}
\end{center}
The point of this table is that as soon as $x$ is moderately large, the approximation $\mathfrak{S}\Li_2(x)$ is exceptionally close to $\pi_2(x)$, and the difference consistently fits the prediction of being of size about $\sqrt{x}$. The twin prime counts are known up to at least $10^{18}$; we have stopped at $10^{15}$ in the interests of space alone!
\subsection{Cram\'er's model}
The Twin Prime Conjecture is supported by various heuristic arguments to justify the main term $\mathfrak{S}\Li_2(x)$. One such heuristic is based on a random model of the primes due to Harold Cram\'er.

If one chooses an integer $n$ uniformly at random between $x$ and $2x$, then the probability that $n$ is prime is roughly $1/\log{x}$ (this is a consequence of the Prime Number Theorem). We might naively guess that the property of $n$ being prime should be roughly independent of the property of $n+2$ being prime when $x$ is large, and so the probability that both $n$ and $n+2$ are prime would be about $1/(\log{n})\cdot 1/\log(n+2)\approx 1/(\log{n})^2$.  This would lead one to guess that there are about $x/(\log{x})^2$ twin primes up to $x$. This is the same order of magnitude as the prediction of Conjecture \ref{cnjctr:TwinPrime}, but misses the constant $\mathfrak{S}$. 

The reason for this discrepancy is because we also need to account for the fact that primes have no small prime factors - if $n>2$ is prime, then we know that $n+2$ is odd, which makes it slightly more likely to be prime. If both $n$ and $n+2$ are primes of size $x$, then they cannot be a multiple of any small prime $p$. If $n$ and $n+2$ were independent, the chances of neither being a multiple of $p$ would be about $(1-1/p)^2$. However, in reality they are not independent; neither is a multiple of $p$ with probability $(1-2/p)$ if $p>2$, since $n$ needs to avoid the residue classes $0$ and $-2 \Mod{p}$. Therefore we should refine the independence assumption, to take this into account. (If $p=2$ then $n$ just needs to avoid the residue class $0\Mod{2}$, which occurs with probability about $1/2$.) This suggests we should have a correction factor of $(1-2/p)/(1-1/p)^2$ for each such small prime $p$, and so the total correction factor should be
\[
\frac{1/2}{(1-1/2)^2}\prod_{p>2}\frac{1-2/p}{(1-1/p)^2}=\mathfrak{S}.
\]
This leads to the guess of $\mathfrak{S}x/\log^2{x}$ twin primes up to $x$, as in Conjecture \ref{cnjctr:TwinPrime}.
\subsection{Circle method heuristic}
An alternative heuristic comes from a technique known as the Hardy-Littlewood circle method. One way to try to prove results like the Twin Prime Conjecture is to try to use the circle method to count the number of twin primes up to $x$. The idea is to first express the count as an integral:
\[
\pi_2(x)=\int_0^1 \Bigl|\sum_{p\le x+2}e^{2\pi i p\theta}\Bigr|^2 e^{4\pi i\theta}d\theta,
\]
and then to try to understand the integral by separating parts when the integrand is large from the parts when it is small. This is the way we prove the three prime version of Goldbach's conjecture, for example. We know that the integrand is large whenever $\theta$ is close to a rational with small denominator. For these parts of the integral, we can understand the relevant exponential sums using results on primes in arithmetic progressions. This allows us to get an asymptotic estimate for the contribution from these $\theta$, which turns out to be
\[
(1+o(1))\mathfrak{S}\int_2^x\frac{dt}{\log^2{t}}.
\]
This is exactly the expected main term! We would successfully prove the Twin Prime Conjecture if we could show the contribution amongst $\theta$ \emph{not} close to a rational with small denominator were small. Unfortunately, we do not know how to do this, although we expect this to be the case. In similar simpler problems we can show these terms are negligible, but for the Twin Prime Problem we don't know how to show cancellation between different $\theta$, and considering each $\theta$ individually is insufficient. However, since we still expect the contribution to be small, this gives an alternative justification for the guess of Conjecture \ref{cnjctr:TwinPrime}.
\section{Sieve methods and related problems}
Although the Twin Prime Conjecture is still very much open, we have succeeded in making partial progress in various directions. The aim of this article is to outline some recent breakthroughs making dramatic new progress on questions related to the conjecture.

\subsection{Early results from sieve methods}
During the $20^{th}$ century, a key development in Prime Number Theory was the invention of \emph{sieve methods}. These are flexible tools for studying the primes, by relating the primes to simpler mathematical objects which can be analyzed effectively. These simpler objects are weighted integers, where the weights are chosen carefully so as to be heavily concentrated on the primes, but are simple enough to calculate with.

Viggo Brun, the pioneer of sieve methods, showed that the number of twin primes up to $x$ is $O(x(\log\log{x}/\log{x})^2)$, which is only slightly larger than the conjectured truth. Slight refinements allow us to show that Conjecture \ref{cnjctr:TwinPrime} cannot underestimate the number of twin primes by more than a bounded amount.
\begin{thrm}[Sieve bound on Twin Primes]
There is a constant $C>0$ such that for $x\ge 3$
\[
\pi_2(x)\le C \mathfrak{S}\int_2^x\frac{dt}{(\log{t})^2}.
\]
\end{thrm}
The current record is that one can take $C$ slightly below $3.4$ if $x$ is large enough (\cite{Wu}, building on the earlier work \cite{Chen2,FouvryIwaniec2,Fouvry2,BFIi,FouvryGrupp}). Even Brun's original result implied that the sum of reciprocals of twin primes converges, unlike the sum of reciprocals of the primes.
\begin{crllry}[Brun's Constant]
The sum of reciprocals of all twin prime pairs converges to a constant
\[
B:=\Bigl(\frac{1}{3}+\frac{1}{5}\Bigr)+\Bigl(\frac{1}{5}+\frac{1}{7}\Bigr)+\Bigl(\frac{1}{11}+\frac{1}{13}\Bigr)+\dots \approx 1.902\dots
\]
\end{crllry}
As well as proving upper bounds for counts of primes with certain primes, sieve methods are also generally able to show one can solve many such problems if one replaces `prime' with `an integer with at most $r$ of prime factors', for some fixed constant $r$. In the context of the Twin Prime Conjecture, the strongest result of this type is Chen's celebrated theorem \cite{Chen}.
\begin{thrm}[Chen's Theorem]
There are infinitely many primes $p$ such that $p+2$ has at most $2$ prime factors.
\end{thrm}
Unfortunately there is a key limitation inherent in standard sieve methods which means very significant new ideas are required if one wants to prove the Twin Prime Conjecture itself. Indeed, it is known that sieve methods struggle to distinguish between integers with an even or an odd number of prime factors, and without overcoming this limitation one can't hope to prove a stronger result than Chen's theorem.

\subsection{Bounded gaps between primes}
An alternative approximation to the Twin Prime Conjecture is to try to show that there are pairs of primes which are close to each other, even if they do not exactly differ by 2. This problem was studied by many authors \cite{Erdos, Rankin, Ricci,BombieriDavenport,Piltai,Uchiyama,Huxley,Huxley2,Huxley3,Maier}, but progress was slow. In 2005 a breakthrough due to Goldston, Pintz and Y\i ld\i r\i m \cite{GPY} showed that there are pairs of primes arbitrarily close to each other compared with the average gap.
\begin{thrm}\label{thrm:GPY}
Let $\epsilon>0$. Then there are infinitely many pairs $(p_1,p_2)$ of distinct primes such that
\[
|p_1-p_2|\le \epsilon \log{p_1}
\]
\end{thrm}
Moreover, the new method - now known as the `GPY method' - could show that there are infinitely many such pairs $(p_1,p_2)$ with $|p_1-p_2|$ bounded, if one assumed a widely believed conjecture on the distribution on primes in arithmetic progressions. Despite many significant results in this direction \cite{FouvryIwaniec,BFIi,BFIii,BFIiii} this conjecture seemed out of reach. In 2013, Zhang \cite{Zhang} proved a suitable variant of this conjecture, which allowed him to use the GPY method to prove that there are infinitely many bounded gaps between primes.
\begin{thrm}[Zhang's Theorem]\label{thrm:Zhang}
There are infinitely many pairs $(p_1,p_2)$ of distinct primes such that
\[
|p_1-p_2|\le \num{70000000}.
\]
\end{thrm}
After Zhang's breakthrough, a new variant of the GPY method was discovered independently by the author \cite{SmallGaps} and Tao (unpublished). This gave an alternative way of proving bounded gaps between primes, but had several other consequences as well since it was flexible and could show the existence of clumps of many primes in bounded length intervals. Moreover, it turned out this new method produced stronger numerical results for the size of the gap. A joint collaborative project was set up to optimize the numerics in this bound. First Zhang's arguments were optimized and improved in \cite{Polymath8a}, and then the new method was also optimized, resulting in the current world record \cite{Polymath}.
\begin{thrm}\label{thrm:Polymath}
There are infinitely many pairs $(p_1,p_2)$ of distinct primes such that
\[
|p_1-p_2|\le 246.
\]
\end{thrm}
All of these approaches actually prove a weak approximation of a far-reaching generalization of the Twin Prime Conjecture, and then deduce results about gaps between primes from this result. The generalization is instead of asking for simultaneous prime values of the functions $n$ and $n+2$, it considers a more general set of many linear functions.

\begin{dfntn}[Admissibility]
Let $\mathcal{L}=\{L_1,\dots,L_k\}$ be a set of linear functions in one variable with integer coefficients and positive lead coefficient (i.e. $L_i(n)=a_in+b_i$ for some $a_i,b_i\in\mathbb{Z}$ with $a_i> 0$). 

We say $\mathcal{L}$ is \emph{admissible} if $\prod_{i=1}^kL_i(n)$ has no fixed prime divisor. (i.e. if for every prime $p$ there is an integer $n_p$ such that $\prod_{i=1}^kL_i(n_p)$ is coprime to $p$.)
\end{dfntn}
Given a set of linear functions, it is easy to check whether they are admissible or not, so this should be thought of as a simple property. If a set of linear functions is not admissible (i.e. there is a fixed prime divisor of $\prod_{i=1}^kL_i(n)$) then it is straightforward to see that there can only be finitely many integers $n$ such that all the $L_i(n)$ are prime. The far-reaching Prime $k$-tuples Conjecture claims that this simple necessary condition is in fact sufficient.

\begin{cnjctr}[Prime $k$-tuples Conjecture]\label{cnjctr:KTuples}
Let $\mathcal{L}=\{L_1,\dots,L_k\}$ be an admissible set of linear functions. Then there are infinitely many integers $n$ such that all of $L_1(n),\dots,L_k(n)$ are primes.
\end{cnjctr}
Conjecture \ref{cnjctr:Twin} follows immediately from this by taking $\mathcal{L}=\{L_1,L_2\}$ with $L_1(n)=n$, $L_2(n)=n+2$, but the Prime $k$-tuples Conjecture goes much further and would have a huge number of applications. It almost completely describes the `small scale' structure of the primes. Since we do not know how to solve the Twin Prime Conjecture, we certainly do not know how to solve the Prime $k$-tuples Conjecture. Recent developments, however, have been based on establishing weak versions of the Prime $k$-tuples Conjecture, where \emph{several} rather than \emph{all} of the linear functions are simultaneously prime. One consequence is we now can prove the following Theorem \cite{Dense}.

\begin{thrm}[Weak $k$-tuples]\label{thrm:WeakTuples}
There is a constant $c>0$ such that the following holds.

Let $\mathcal{L}$ be an admissible set of $k$ linear functions. Then there are infinitely many integers $n$ such that at least $c\log{k}$ of $L_1(n),\dots,L_k(n)$ are primes.
\end{thrm}
The fact that there are infinitely many bounded gaps between primes follows immediately from Theorem \ref{thrm:WeakTuples} on choosing $k$ such that $c\log{k}\ge 2$ and $\mathcal{L}=\{L_1,\dots,L_k\}$ with $L_i(n)=n+h_i$ for suitable constants $h_i$ such that $\mathcal{L}$ is admissible (choosing $h_i=i\times k!$ works, for example). The stronger numerical results such as  Theorem \ref{thrm:Polymath} come from numerically optimizing the proof of Theorem \ref{thrm:WeakTuples} to determine what the smallest value of $k$ is such that $c\log{k}$ can be replaced by $2$.

We briefly overview some of the ideas behind the GPY method in Section \ref{sec:GPY}, and then sketch some of the ideas in the proof of Theorem \ref{thrm:WeakTuples} in Section \ref{sec:Sieve}.

\subsection{Chowla's conjecture}

Even though we can prove results such as Theorem \ref{thrm:Polymath}, significant new ideas are required to prove the Twin Prime Conjecture itself. As mentioned before, sieve methods have a limitation in that they struggle to distinguish  numbers with an even number of prime factors from numbers with an odd number of prime factors. Results such as Theorem \ref{thrm:WeakTuples} exploit the flexibility of only asking for some of the linear functions to be prime to sidestep this limitation, without directly addressing it. To prove the Twin Prime Conjecture, however, it is likely one will need to address this limitation directly.

One way of seeing this is to consider variants of questions about the primes where instead one considers the Liouville $\lambda$ function, given by
\[
\lambda(n)=\begin{cases}
-1,\qquad &\text{ $n$ has an odd number of prime factors,}\\
1,&\text{otherwise.}
\end{cases}
\]
Showing cancellation in $\lambda(n)$ was generally considered as difficult as studying the primes in any given problem. For example, Chowla \cite{Chowla} made a general conjecture in the same spirit as Conjecture \ref{cnjctr:KTuples}. For the analogue of the Twin Prime Conjecture, this takes the following form.
\begin{cnjctr}[Weak form of Chowla's Conjecture]\label{cnjctr:Chowla}
We have
\[
\sum_{n\le x}\lambda(n)\lambda(n+2)=o(x).
\]
\end{cnjctr}
Until very recently it was believed that this conjecture was essentially as hard as the Twin Prime Conjecture, and the limitations of sieve methods meant that they were poorly suited to addressing this problem. It was therefore considered a major advance when Tao \cite{Tao} showed that a slightly weaker version of Conjecture \ref{cnjctr:Chowla} is true, where we use a logarithmic average instead of a standard average.
\begin{thrm}[Logarithmically Averaged Chowla]\label{thrm:LogChowla}
\[
\sum_{n\le x}\frac{\lambda(n)\lambda(n+2)}{n}=o(\log{x}).
\]
\end{thrm}
This result introduced several new ideas, but also relied heavily on a recent breakthrough of Matom\"aki and Radziwi\l\l\,\cite{MatomakiRadziwill}, who had shown that for the question of almost all short intervals, one could prove results for $\lambda(n)$ which are much better than what is known for primes.
\begin{thrm}[Liouville in short intervals]\label{thrm:MatomakiRadziwill}
For almost all intervals $[x,x+h]\subseteq[Y,2Y]$, we have
\[
\sum_{n\in[x,x+h]}\lambda(n)=O\Bigl(\frac{h}{(\log{h})^{1/10}}\Bigr).
\]
\end{thrm}
The Matom\"aki-Radziwi\l\l\,Theorem is optimal in the sense that provided $h\rightarrow \infty$ there is cancellation in the Liouville function in almost all short intervals. By way of comparison, for primes we only know that there is the expected number of primes in almost all intervals $[x,x+h]\subseteq[Y,2Y]$ when $h\ge Y^{1/6}$, a much weaker result. This therefore demonstrates that the previous belief that many questions about $\lambda(n)$ were as hard as the same question about the primes is not correct.

Unfortunately it currently appears that these breakthroughs do not directly address the Twin Prime Conjecture itself - they rely on showing that the behaviour of small prime factors and large prime factors is roughly independent, and cancellation comes from small prime factors of integers. This means the methods do not appear to adapt to say anything about the primes (where there are no small prime factors), but the mere fact that something can be said about these questions which directly address the limitations of sieve techniques is a big advance.

Although new ideas appear to be required to prove the Twin Prime Conjecture, the breakthroughs of Theorems \ref{thrm:WeakTuples}, \ref{thrm:LogChowla} and \ref{thrm:MatomakiRadziwill} have already had several further applications, including new results on large gaps between primes \cite{FGKT,FGKMT,MaynardLarge}, the resolution of the Erd\H os discrepancy problem \cite{ErdosDiscrepency}, as well as many other results on the distribution of primes \cite{Joni2,Dense,Thorner,LemkeOliver,Freiberg,FreibergII,Pollack,Hongze,BakerPollack,MatomakiShao,BakerZhao,ChuaParkSmith,Vatwani,Troupe,PintzII,PintzIII,Huang,BakerZhaoII,BanksFreiberg,BakerFreiberg,Kaptan,Parshall,PollackThompson,BakerFreiberg,MaierRassias,Pratt} and correlations of multiplicative functions \cite{Joni,TaoKaisaMaks,Klurman1,Klurman2,Goudout,LesterMR,FranziHost,Franzi}.

We will sketch some of the key ideas behind Theorem \ref{thrm:LogChowla} in Section \ref{sec:Tao}, assuming a variant of Theorem \ref{thrm:MatomakiRadziwill}. We then sketch some of the ideas behind Theorem \ref{thrm:MatomakiRadziwill} in Section \ref{sec:MatomakiRadziwill}.
\section{The GPY method and Weak Prime \texorpdfstring{$k$}{k}-tuples}\label{sec:GPY}
We aim to prove Theorem \ref{thrm:WeakTuples} by the `GPY method', which can be interpreted as a first moment method which was introduced to study small gaps between primes by Goldston, Pintz and Y\i ld\i r\i m \cite{GPY}. This takes the following basic steps, for some given set $\{L_1,\dots,L_k\}$ of distinct functions satisfying the hypotheses of Conjecture \ref{cnjctr:KTuples}:

\begin{enumerate}
\item[1.] We choose a probability measure $w$ supported on integers in $[X,2X]$.
\item[2.] We calculate the expected number of the functions $L_i(n)$ which are prime at $n$, if $n$ is chosen randomly with probability $w(n)$.
\item[3.] If this expectation is at least $m$, then there must be some $n\in[X,2X]$ such that at least $m$ of $L_1(n),\dots, L_k(n)$ are primes.
\item[4.] If this holds for all large $X$, then there are infinitely many such integers $n$.
\end{enumerate}
This procedure only works if we can find a probability measure $w$ which is suitably concentrated on integers $n$ where many of the $L_i(n)$ are prime, but at the same time is simple enough that we can calculate this expectation unconditionally. We note that by linearity of expectation, it suffices to be able to calculate the probability that any one of the linear functions is prime. However, any slowly changing smooth function $w$ is insufficient since the primes have density 0 in the integers, whereas any choice of $w$ explicitly depending on the joint distribution of prime values of the $L_i$ is likely to be too complicated to handle unconditionally. 

Sieve methods provide natural choices for the probability measure $w$, since they produce arithmetic objects concentrated on primes but simple enough to calculate with. The situation of simultaneous prime values of $L_1,\dots,L_k$ is a `$k$-dimensional' sieve problem. For such problems when $k$ is large but fixed, the Selberg sieve tends to be the type of sieve which performs best. The standard choice of Selberg sieve weights (which are essentially optimal in closely related situations) are
\begin{equation}
w(n)\propto\Bigl(\sum_{\substack{d|\prod_{i=1}^k L_i(n) \\ d<R}}\mu(d)\Bigl(\log\frac{R}{d}\Bigr)^k\Bigr)^2,\label{eq:SelbergWeights}
\end{equation}
where $R$ is a parameter which controls the complexity of the sieve weights, and $w$ is normalized to sum to $1$ on $[X,2X]$.

To calculate the probability that $L_j(n)$ is prime with this choice of $w(n)$, we wish to estimate the sum of $w(n)$ over $n\in[X,2X]$ such that $L_j(n)$ is prime. To do this we typically expand the divisor sum in the definition \eqref{eq:SelbergWeights} and swap the order of summation. This reduces the problem to estimating the number of prime values of $L_j(n)$ for $n\in [X,2X]$ in many different arithmetic progressions with moduli of size about $R^2$. The Elliott-Halberstam conjecture \cite{ElliottHalberstam} asserts that we should be able to do this when $R^2<X^{1-\epsilon}$, but unconditionally we only know how to do this when $R^2<X^{1/2-\epsilon}$, using the Bombieri-Vinogradov Theorem \cite{Bombieri,Vinogradov}. After some computation one finds that, provided we do have suitable estimates for primes in arithmetic progressions, the choice \eqref{eq:SelbergWeights} gives
\begin{equation}
\mathbb{E}\,\#\{1\le i\le k:L_i(n)\text{ prime}\}=\Bigl(2-\frac{2}{k+1}+o(1)\Bigr)\frac{\log{R}}{\log{X}}. \label{eq:SelbergExpectation}
\end{equation}
In particular, this is less than $1$ for all large $X$, even if we optimistically assume the Elliott-Halberstam conjecture and take $R\approx X^{1/2-\epsilon}$, and so it appears that we will not be able to conclude anything about primes in this manner.

The groundbreaking work of Goldston-Pintz-Y\i ld\i r\i m \cite{GPY} showed that a variant of this choice of weight actually performs much better. They considered
\begin{equation}
w(n)\propto\Bigl(\sum_{\substack{d|\prod_{i=1}^k L_i(n) \\ d<R}}\mu(d)\Bigl(\log\frac{R}{d}\Bigr)^{k+\ell}\Bigr)^2,\label{eq:GPYWeights}
\end{equation}
where $\ell$ is an additional parameter to be optimized over. With the choice \eqref{eq:GPYWeights}, one finds that provided we have suitable estimates for primes in arithmetic progressions, we obtain
\[
\mathbb{E}\,\#\{1\le i\le k:L_i(n)\text{ prime}\}=\Bigl(4+O\Bigl(\frac{1}{\ell}\Bigr)+O\Bigl(\frac{\ell}{k}\Bigr)\Bigr)\frac{\log{R}}{\log{X}}.
\]
This improves upon \eqref{eq:SelbergExpectation} by a factor of about 2 when $k$ is large and $\ell\approx k^{1/2}$. This falls just short of proving that two of the linear functions are simultaneously prime when $R=X^{1/4-\epsilon}$ as allowed by the Bombieri-Vinogradov theorem, but any small improvement allowing $R=X^{1/4+\epsilon}$ would give bounded gaps between primes! By considering additional possible primes, Goldston Pintz and Y\i ld\i r\i m were to use this argument to prove Theorem \ref{thrm:GPY}.
 
Building on earlier work, Zhang \cite{Zhang} succeeded in establishing an extended version of the Bombieri-Vinogradov Theorem for moduli with no large prime factors allowing $R=X^{1/4+\epsilon}$, and ultimately this allowed him to show
\[
\mathbb{E}\,\#\{1\le i\le k:L_i(n)\text{ prime}\}>1,
\]
if $k>\num{3500000}$ and $X$ is sufficiently large. The key breakthrough in Zhang's work was this result on primes in arithmetic progressions of modulus slightly larger than $X^{1/2}$. By choosing suitable linear functions, this then established Theorem \ref{thrm:Zhang}.

An alternative approach to extending the work of Goldston, Pintz and Y\i ld\i r\i m  was developed independently by the author \cite{SmallGaps} and Tao (unpublished). The key difference was to consider the multidimensional generalization%\footnote{Selberg had Goldston and Y\i ld\i r\i m also had earlier work which considered similar setups, which were not quite sufficient for bounded gaps between primes}
\begin{equation}
w(n)\propto\Bigl(\sum_{\substack{d_1,\dots,d_k \\ d_i|L_i(n) \\ \prod_{i=1}^k d_i<R}}\mu(d_1\cdots d_k)F\Bigl(\frac{\log{d_1}}{\log{R}},\dots,\frac{\log{d_k}}{\log{R}}\Bigr)\Bigr)^2,\label{eq:NewWeights}
\end{equation}
for suitable smooth functions $F:\mathbb{R}^k\rightarrow \mathbb{R}$ supported on $[0,\infty)^k$. The flexibility of allowing the function $F$ to depend on each divisor $d_1,\dots,d_k$ of $L_1(n),\dots,L_k(n)$ allows us to make $w(n)$ more concentrated on integers $n$ when many of the $L_i(n)$ are prime (at the cost of making it less concentrated on integers $n$ when \emph{all} the linear functions are prime). This is the key new idea allowing for the proof of Theorem \ref{thrm:WeakTuples}.

For further reading on the GPY method, we recommend the excellent survey article of Soundararajan \cite{Sound}. For further reading about the more recent developments surrounding bounded gaps between primes, we recommend the survey articles of Granville \cite{Granville} and Kowalski \cite{Kowalski}.

\section{A modified GPY sieve method}\label{sec:Sieve}

For a given fixed admissible set $\{L_1,\dots,L_k\}$, we wish to estimate
\[
\mathbb{E}\#\{1\le i\le k:\, L_i(n)\text{ prime}\}=\sum_{i=1}^k\mathbb{P}(L_i(n)\text{ prime}),
\]
when $n\in [x,2x]$ is chosen randomly with probability according to \eqref{eq:NewWeights}. Let
\[
\nu(n)=\Bigl(\sum_{\substack{d_1,\dots,d_k \\ d_i|L_i(n) \\ \prod_{i=1}^k d_i<R}}\mu(d_1\cdots d_k)F\Bigl(\frac{\log{d_1}}{\log{R}},\dots,\frac{\log{d_k}}{\log{R}}\Bigr)\Bigr)^2,
\]
for some piecewise smooth function $F$, independent of $x$. Thus $w(n)=\nu(n)/(\sum_{m\in[x,2x]}\nu(m) )$, and we have
\[
\mathbb{P}(L_i(n)\text{ prime})=\frac{\sum_{n\in[x,2x]}\nu(n)\mathbf{1}_{prime}(L_i(n))}{\sum_{n\in[x,2x]}\nu(n)}.
\]
Therefore it suffices to get asymptotic estimates for the numerator and denominator, and show that this is large when $x$ is large and $F=F_k$ is chosen suitably.

\subsection{First sieve calculations}
We first try to estimate $\sum_{n\in[x,2x]}\nu(n)$. Expanding the square in the definition of $\nu(n)$ and swapping the order of summation, we find that this is given by
\begin{equation}
\sum_{\substack{d_1,\dots,d_k  \\ \prod_{i=1}^k d_i<R}}\sum_{\substack{e_1,\dots,e_k \\ \prod_{i=1}^k e_i<R}}\mu(\mathbf{d})\mu(\mathbf{e})F\Bigl(\frac{\log{\mathbf{d}}}{\log{R}}\Bigr)F\Bigl(\frac{\log{\mathbf{e}}}{\log{R}}\Bigr)\sum_{\substack{n\in[x,2x]\\ d_i|L_i(n)\,\forall i\\ e_i|L_i(n)\forall i}}1,
\label{eq:Expanded}
\end{equation}
where we have written $\mu(\mathbf{d})=\mu(d_1\cdots d_k)$ and 
\[
F\Bigl(\frac{\log{\mathbf{d}}}{\log{R}}\Bigr)=F\Bigl(\frac{\log{d_1}}{\log{R}},\dots,\frac{\log{d_k}}{\log{R}}\Bigr)
\]
for notational simplicity.

We see that the conditions $d_i|L_i(n)$ and $e_i|L_i(n)$ mean that $\lcm(d_i,e_i)|L_i(n)$. Imagine for simplicity that $L_i(n)=a_in+b_i$ are such that $\prod_{i<j}(b_i-b_j)$ divides $a_\ell$ for all $1\le \ell \le k$ (this avoids a few minor issues with small primes). Let $A=\lcm(a_1,\dots,a_k)$. Then, by the Chinese Remainder Theorem,  we can rewrite the conditions
\[
\lcm(d_i,e_i)|L_i(n)\quad \text{for }1\le i\le k
\]
as
\[
n\equiv n_0\Mod{A\prod_{i=1}^k\lcm(d_i,e_i)}
\]
for some $n_0$, provided $A,\lcm(d_1,e_1),\dots,\lcm(d_{k-1},e_{k-1})$ and $\lcm(d_k,e_k)$ are all pairwise coprime. If they are not pairwise coprime, then there are no solutions to the conditions.

Therefore, we find that the inner sum of \eqref{eq:Expanded} simply counts integers in an arithmetic progression whenever it is non-empty. The number of $n\in[x,2x]$ with $n\equiv a\Mod{q}$ is $x/q+O(1)$. Thus \eqref{eq:Expanded} simplifies to
\[
\frac{x}{A}\sideset{}{^*}\sum_{\substack{d_1,\dots,d_k  \\ \prod_{i=1}^k d_i<R}}\,\,\sideset{}{^*}\sum_{\substack{e_1,\dots,e_k \\ \prod_{i=1}^k e_i<R}}\frac{\mu(\mathbf{d})\mu(\mathbf{e})}{\prod_{i=1}^k\lcm(d_i,e_i)}F\Bigl(\frac{\log{\mathbf{d}}}{\log{R}}\Bigr)F\Bigl(\frac{\log{\mathbf{e}}}{\log{R}}\Bigr)+O(R^2).
\]
Here the asterisk indicates we have the restriction that $a_1\dots a_k,d_1e_1,\dots, d_ke_k$ are pairwise coprime.

This is now a smoothed sum of multiplicative functions, which are standard sums to evaluate in analytic number theory, typically via complex analysis.
\begin{lmm}[Estimate of smoothed sum of multiplicative functions]\label{lmm:Summation}
Let
\[
S=\sideset{}{^*}\sum_{\substack{d_1,\dots,d_k  \\ \prod_{i=1}^k d_i<R}}\,\,\sideset{}{^*}\sum_{\substack{e_1,\dots,e_k \\ \prod_{i=1}^k e_i<R}}\frac{\mu(\mathbf{d})\mu(\mathbf{e})}{\prod_{i=1}^k\lcm(d_i,e_i)}F\Bigl(\frac{\log{\mathbf{d}}}{\log{R}}\Bigr)F\Bigl(\frac{\log{\mathbf{e}}}{\log{R}}\Bigr).
\]
Then
\[
S=(1+o(1))\frac{A^k}{\phi(A)^k(\log{R})^k}\int_{\substack{t_1,\dots,t_k\ge 0\\ t_1+\dots +t_k\le 1}}\frac{\partial^k F(t_1,\dots,t_k)}{\partial t_1\dots \partial t_k}^2dt_1\dots dt_k.
\]
\end{lmm}
Therefore, letting $\tilde{I}(F)$ be the integral above, we see that
\begin{equation}
\sum_{n\in[x,2x]}\nu(n)=(1+o(1))\frac{xA^{k-1}}{\phi(A)^k(\log{R})^k}\tilde{I}(F).
\label{eq:ISum}
\end{equation}
This gives our asymptotic evaluation of $\sum_n\nu(n)$.
\subsection{Sieve calculations with primes}
We now want to calculate $\sum_{n\in[x,2x]}\nu(n)\mathbf{1}_{\text{prime}}(L_j(n))$ for each $j$. For notational simplicity we will consider the case $j=k$, the other cases being essentially identical. As before, expanding out we find this is given by
\begin{equation}
\sum_{\substack{d_1,\dots,d_k  \\ \prod_{i=1}^k d_i<R}}\sum_{\substack{e_1,\dots,e_k \\  \prod_{i=1}^k e_i<R}}\mu(\mathbf{d})\mu(\mathbf{e})F\Bigl(\frac{\log{\mathbf{d}}}{\log{R}}\Bigr)F\Bigl(\frac{\log{\mathbf{e}}}{\log{R}}\Bigr)\sum_{\substack{n\in[x,2x]\\ d_i|L_i(n)\,\forall i\\ e_i|L_i(n)\forall i}}\mathbf{1}_{\text{prime}}(L_k(n)).
\label{eq:Expanded2}
\end{equation}
Since we have restricted $L_k(n)$ to be prime, we must have $d_k=e_k=1$ since $R<x$. As before, we can combine the remaining conditions $d_i,e_i|L_i(n)$ into a single congruence via the Chinese Remainder Theorem
\[
n\equiv n_0'\Mod{A\prod_{i=1}^{k-1}\lcm(d_i,e_i)},
\]
for some $n_0'$ provided $A$, $d_1e_1$, $\dots$, $d_{k-1}e_{k-1}$ are pairwise coprime (and if they are not the inner sum is empty). Moreover, we can check that $n_0'$ is coprime to the modulus. Therefore we wish to estimate primes in arithmetic progressions. Dirichlet's famous theorem (and its refinement by Siegel and Walfisz) allows us to count primes in arithmetic progressions, but only when the modulus is very small. Unfortunately here the modulus can be as large as $R^2$, which we want to be a power of $x$, which is much too large to handle without making significant progress on the issue of Siegel zeros. Fortunately, we only need an average form of this result, and the celebrated Bombieri-Vinogradov Theorem \cite{Bombieri,Vinogradov} provides an adequate substitute for the Generalized Riemann Hypothesis here.

\begin{thrm}[Bombieri-Vinogradov Theorem]
Let $\pi(x,q,a)$ denote the number of primes $p\le x$ with $p\equiv a\Mod{q}$. Then for any $A,\epsilon>0$ there is a constant $C(A,\epsilon)>0$ such that
\[
\sum_{q\le x^{1/2-\epsilon}}\sup_{(a,q)=1}\Bigl|\pi(x,q,a)-\frac{\pi(x)}{\phi(q)}\Bigr|\le \frac{C(A,\epsilon) x}{(\log{x})^A}.
\]
\end{thrm}
This gives a good result provided all our moduli are bounded by $x^{1/2-\epsilon}$, which means we require $R^2\le x^{1/2-\epsilon}$. Using this, we find that \eqref{eq:Expanded2} is essentially given by
\[
\frac{\pi(x)}{\phi(A)}\sideset{}{^*}\sum_{\substack{d_1,\dots,d_{k-1}  \\ \prod_{i=1}^k d_i<R \\ d_k=1}}\,\,\sideset{}{^*}\sum_{\substack{e_1,\dots,e_{k-1} \\ \prod_{i=1}^k e_i<R \\ e_k=1}}\frac{\mu(\mathbf{d})\mu(\mathbf{e})}{\prod_{i=1}^{k-1}\phi(\lcm(d_i,e_i))}F\Bigl(\frac{\log{\mathbf{d}}}{\log{R}}\Bigr)F\Bigl(\frac{\log{\mathbf{e}}}{\log{R}}\Bigr),
\]
where the asterisks indicate the summation is subject to $A$, $d_1e_1$, $d_2e_2$, $\dots$, $d_{k-1}e_{k-1}$ being coprime. For the purposes of the summation, there is essentially nothing lost if we replace $\phi(\lcm(d_i,e_i))$ with $\lcm(d_i,e_i)$. Thus, a variant of Lemma \ref{lmm:Summation} gives
\begin{align*}
&\sideset{}{^*}\sum_{\substack{d_1,\dots,d_{k-1}  \\ \prod_{i=1}^k d_i<R \\ d_k=1}}\,\,\sideset{}{^*}\sum_{\substack{e_1,\dots,e_{k-1} \\ \prod_{i=1}^k e_i<R\\ e_k=1}}\frac{\mu(\mathbf{d})\mu(\mathbf{e})}{\prod_{i=1}^{k-1}\phi(\lcm(d_i,e_i))}F\Bigl(\frac{\log{\mathbf{d}}}{\log{R}}\Bigr)F\Bigl(\frac{\log{\mathbf{e}}}{\log{R}}\Bigr)\\
&=\frac{(1+o(1))A^{k-1}}{\phi(A)^{k-1}(\log{R})^{k-1}}\int_{\substack{t_1,\dots,t_{k-1}\ge 0\\ t_1+\dots +t_{k-1}\le 1}}\frac{\partial^{k-1} F(t_1,\dots,t_{k-1},0)}{\partial t_1\dots \partial t_{k-1}}^2dt_1\dots dt_{k-1}.
\end{align*}
Letting $\tilde{J}_k(F)$ denote the integral above, and using the Prime Number Theorem $\pi(x)=(1+o(1))x/\log{x}$, we see this gives
\begin{align}
\sum_{n\in[x,2x]}\nu(n)\mathbf{1}_{\text{prime}}(L_k(n))=\frac{(1+o(1))A^{k-1} x}{\phi(A)^k (\log{R})^{k-1}\log{x}}\tilde{J}_k(F).\label{eq:JSum}
\end{align}
Thus, combining \eqref{eq:ISum} and \eqref{eq:JSum}, we find that 
\begin{align*}
\mathbb{P}(L_k(n)\text{ prime})&=\frac{\sum_{n\in[x,2x]}\nu(n)\mathbf{1}_{\text{prime}}(L_k(n))}{\sum_{n\in[x,2x]}\nu(n)}\\
&=(1+o(1))\frac{\log{R}}{\log{x}}\frac{\tilde{J}_k(F)}{\tilde{I}(F)}.
\end{align*}
\subsection{Choice of \texorpdfstring{$F$}{F}}
It will be convenient to work with
\[
\tilde{F}=\frac{\partial^k}{\partial t_1\dots \partial t_k}F.
\]
With this transform, the previous calculations show that for $R\le x^{1/4-\epsilon}$
\begin{equation}
\mathbb{E}\,\#\{1\le i\le k:L_i(n)\text{ prime}\}=\Bigl(\frac{\sum_{i=1}^k J_i(\tilde{F})}{I(\tilde{F})}+o(1)\Bigr)\frac{\log{R}}{\log{X}},
\label{eq:NewExpectation}
\end{equation}
where the $J_i(\tilde{F})$ and $I(\tilde{F})$ are given by
\begin{align*}
J_\ell(\tilde{F}) &=\idotsint\limits_{\sum_{i\ne \ell} t_i\le 1}\Bigl(\int_0^{1-\sum_{i\ne \ell}t_i} \tilde{F}(t_1,\dots,t_k)dt_\ell \Bigr)^2 dt_1\dots dt_{\ell-1}dt_{\ell+1}\dots dt_k,\\
I(\tilde{F}) &=\idotsint\limits_{\sum_{i=1}^k t_i\le 1}\tilde{F}(t_1,\dots,t_k)^2dt_1\dots dt_k.
\end{align*}
In particular, we can show that many of the $L_i$ are simultaneously prime infinitely often, if we can show that $\sup_{\tilde{F}}(\sum_{i=1}^k J_i(\tilde{F})/I(\tilde{F}))\rightarrow \infty$ as $k\rightarrow\infty$, where the supremum is over all piecewise smooth $\tilde{F}$.

A key advantage of this generalization is that we can make use of high dimensional phenomena such as concentration of measure. We concentrate on functions $\tilde{F}$ the form
\[
\tilde{F}(t_1,\dots,t_k)=\begin{cases}
\prod_{i=1}^k G(t_i),\qquad &\sum_{i=1}^kt_i<1,\\
0,&\text{otherwise,}
\end{cases}
\]
for some smooth function $G(t)$ (depending on $k$). This allows us to make a probabilistic interpretation of the integrals $J_i(\tilde{F})$ and $I(\tilde{F})$. Let $Z_1,\dots,Z_k$ be i.i.d. random variables (depending on $k$) on $[0,\infty]$ with probability density function $G^2$, expectation $\mu=\int_0^\infty t G(t)^2dt$ and variance $\sigma^2=\int_0^\infty (t-\mu)^2G(t)^2dt$. Then
\begin{align*}
I(\tilde{F})&=\mathbb{P}\Bigl(\sum_{i=1}^k Z_i<1\Bigr),\\
J_\ell(\tilde{F})&\ge \Bigl(\int_0^{1/2}G(t)dt\Bigr)^2\mathbb{P}\Bigl(\sum_{i=1}^k Z_i<\frac{1}{2}\Bigr).
\end{align*}
The random variable $\sum_{i=1}^k Z_i$ has mean $k\mu$ and variance $k\sigma^2$, and so it becomes concentrated on $k\mu$ when $k$ is large provided $\sigma/k^{1/2}\mu\rightarrow 0$ as $k\rightarrow\infty$. In particular, if $\mu< 1/3k$ and $\sigma^2k\rightarrow 0$, then $\mathbb{P}(\sum_{i=1}^k Z_i<1/2)$ approaches 1 as $k\rightarrow\infty$. Therefore, to show that $\sum_{\ell=1}^kJ_\ell(\tilde{F})/I(\tilde{F})\rightarrow\infty$ as $k\rightarrow\infty$, it suffices to find a function $G$ satisfying:
\begin{align*}
\int_0^\infty t G(t)^2&<\frac{1}{3k},\\
\int_0^\infty G(t)^2 dt&=1,\\
k\int_0^\infty t^2G(t)^2dt&\rightarrow 0,\\
k\Bigl(\int_0^\infty G(t)dt\Bigr)^2&\rightarrow \infty \quad \text{as}\quad k\rightarrow \infty.
\end{align*}
We find that choosing
\begin{equation}
G(t)\approx \begin{cases}
\frac{\sqrt{k\log{k}}}{1+t k\log{k}},\qquad &t<k^{-3/4},\\
0,&\text{otherwise,}
\end{cases}\label{eq:SmoothChoice}
\end{equation}
gives a function $G$ satisfying these constraints. (This choice can be found via the calculus of variations, and further calculations show that this choice is essentially optimal.) Putting this all together, we find that for some constant $c>0$
\[
\mathbb{E}\,\#\{i:\,L_i(n)\text{ prime}\}\ge ( c\log{k}+o(1))\frac{\log{R}}{\log{X}}.
\]
In particular, taking $R=X^{1/4-o(1)}$ (as allowed by the Bombieri-Vinogradov theorem) and letting $k$ be sufficiently large, we find that there are infinitely many integers $n$ such that $(c/4+o(1))\log{k}$ of the $L_i(n)$ are simultaneously prime. Performing these calculations carefully allows one to take $c\approx 1$ when $k$ is large.

Morally, the effect of such a choice of function $\tilde{F}$ is to make `typical' divisors $(d_1,\dots,d_k)$ occurring in \eqref{eq:NewWeights} to have $\prod_{i=1}^k d_i$ smaller, but for it to be more common that some of the components of $(d_1,\dots,d_k)$ are unusually large when compared with \eqref{eq:SelbergWeights} or \eqref{eq:GPYWeights}. This correspondingly causes the random integer $n$ chosen with probability $w(n)$ to be such that the $L_i(n)$ are more likely to have slightly smaller prime factors, but it is also more likely that some of the $L_i(n)$ are prime.

\subsection{Consequences for small gaps between primes}

The above argument shows that there is a constant $c>0$ such that for any set $\{L_1,\dots,L_k\}$ of integral linear functions satisfying the hypotheses of Conjecture \ref{cnjctr:KTuples}, at least $c\log{k}$ of the linear functions $L_i$ are simultaneously prime infinitely often. To show that there are primes close together, we simply take the linear functions to be of the form $L_i(n)=n+h_i$ for some integers $h_1< \dots < h_k$ chosen to make $\prod_{i=1}^k L_i(n)$ have no fixed prime divisor, and so that $h_k-h_1$ is small. 

In general, a good choice of the $h_i$ is to take $h_i$ to be the $i^{th}$ prime after the integer $k$. With this choice, $\prod_{i=1}^kL_i(n)$ has no fixed prime divisor and $h_k-h_1\approx k\log{k}$, so we can find intervals of length $k\log{k}$ containing $c\log{k}$ primes infinitely often. By working out the best possible implied constants, the argument we have sketched allows us to show that there are $m$ primes in an interval of length $O(m^3e^{4m})$ infinitely often. The fact that we can find many primes rather than just 2, and we can do so taking $R=x^{\epsilon}$ has important other applications, such as for large gaps between primes \cite{MaynardLarge,FGKT,FGKMT}.

\subsection{Numerical optimization}
If we are only interested in just how small a \textit{single} gap can be, then we can improve the analysis and get an explicit bound on the size of the gap by adopting a numerical analysis perspective. The key issue is to find the smallest value of $k$ such that for all large $X$
\[
\mathbb{E}\,\#\{1\le i\le k:\,L_i(n)\text{ prime}\}>1,
\]
since this immediately implies that two of the linear functions are simultaneously prime. Recalling that we can choose $R=x^{1/4-\epsilon}$ for any $\epsilon>0$, using \eqref{eq:NewExpectation} we reduce the problem to finding a value of $k$ as small as possible, such that we can find a function $\tilde{F}:[0,\infty)^k\rightarrow \mathbb{R}$ with
\[
\frac{\sum_{i=1}^k J_i(\tilde{F})}{I(\tilde{F})}>4.
\]
We fix some basis functions $g_1,\dots,g_r:[0,\infty)^k\rightarrow \mathbb{R}$ which are supported on $(t_1,\dots,t_k)$ satisfying $\sum_{i=1}^k t_i\le 1$, and restrict our attention to functions $\tilde{F}$ in the linear span of the $g_i$; that is functions $\tilde{F}$ of the form
\[
\tilde{F}(t_1,\dots,t_k)=\sum_{i=1}^r f_i g_i(t_1,\dots,t_k),
\]
for some coefficients $\mathbf{f}=(f_1,\dots f_r)\in\mathbb{R}^r$ which we think of as variables we will optimize over. For such a choice, we find that $\sum_{i=1}^k J_i(\tilde{F})$ and $I(\tilde{F})$ are both quadratic forms in the variables $f_1,\dots,f_r$, and the coefficients of these quadratic forms are explicit integrals in terms of $g_1,\dots,g_r$. If we choose a suitably nice basis $g_1,\dots,g_r$, then these integrals are explicitly computable, and so we obtain explicit $r\times r$ real symmetric matrices $M_1$ and $M_2$ such that
\[
\frac{\sum_{i=1}^k J_i(\tilde{F})}{I(\tilde{F})}=\frac{\mathbf{f}^T M_1\mathbf{f}}{\mathbf{f}^T M_2\mathbf{f}}.
\]
We then find that the choice of coefficients $\mathbf{f}$ which maximizes this ratio is the eigenvector of $M_2^{-1}M_1$ corresponding to the largest eigenvalue, and the value of the ratio is given by this largest eigenvalue. Thus the existence of a good function $\tilde{F}$ can be reduced to checking whether the largest eigenvalue of a finite matrix is larger than 4, which can be performed numerically by a computer.

  If we choose the $g_i$ to be symmetric polynomials of low degree, then this provides a nice basis since the corresponding integrals have a closed form solution, and allows one to make arbitrarily accurate numerical approximations to the optimal function $\tilde{F}$ with enough computation. For the problem at hand, these numerical calculations are large but feasible. This approach ultimately allows one to show that if $k=54$ there is a function $\tilde{F}$ such that $\sum_{i=1}^{54}J_i(\tilde{F})/I(\tilde{F})>4$. 
  
  To turn this into small gaps between primes we need to choose the shifts $h_i$ in our linear functions $L_i(n)=n+h_i$ so that $\prod_{i=1}^{54}L_i(n)$ has no fixed prime divisor, and the $h_i$ are in as short an interval as possible. This is a feasible exhaustive numerical optimization problem, with an optimal choice of the $\{h_1,\dots,h_{54}\}$ given by \footnote{Such computations were first performed by Engelsma - see \url{http://www.opertech.com/primes/k-tuples.html}}
  \begin{align*}
  \{&0, 2, 6, 12, 20, 26, 30, 32, 42, 56, 60, 62, 72, 74, 84, 86, 90, 96, \\
  &\quad  104, 110, 114, 116, 120, 126, 132, 134, 140, 144, 152, 156, \\
  &\quad 162, 170, 174, 176, 182, 186, 194, 200, 204, 210, 216, 222, \\
  &\quad  224, 230, 236, 240, 242, 246, 252, 254, 260, 264, 266, 270\}.
  \end{align*}
  Putting this together, we find that there are infinitely many pairs of primes $p\ne q$ such that
\[
|p-q|\le 270.
\]
This result is not quite the current record - in \cite{Polymath} we make some further technical refinements to the sieve, which corresponds to modifying the expressions $J_i(\tilde{F})$ and $I(\tilde{F})$ slightly. This ultimately allows us to improve $k=54$ to $k=50$ and correspondingly improve from gaps of size at most 270 to gaps of at most 246, giving Theorem \ref{thrm:Polymath}. The main ideas are the same as above.

\section{Logarithmic Chowla Conjecture}\label{sec:Tao}
We wish to sketch the proof of the following result due to Tao. For minor technical reasons, we will consider the following variant of Theorem \ref{thrm:LogChowla}, where we replace $n+2$ with $n+1$:
\begin{thrm}[Variant of logarithmically averaged 2-point Chowla conjecture]
As $x\rightarrow \infty$, we have
\[
\sum_{n<x}\frac{\lambda(n)\lambda(n+1)}{n}=o(\log{x})
\]
\end{thrm}
This will rely heavily on the Matom\"aki-Radziwi\l\l\,Theorem (Theorem \ref{thrm:MatomakiRadziwill}), which we talk about in the next section. In fact, for simplicity we will use a slight variant of their theorem, which is the following, established in \cite{MatomakiRadziwillTao}.
\begin{thrm}[Extended Matom\"aki-Radziwi\l\l]\label{thrm:Extended}
Let $\theta\in\mathbb{R}$. Then for all but at most $X/(\log{h})^{1/10}$ integers $n\le X$ we have
\[
\Bigl|\sum_{j\le h}\lambda(n+j)e(j\theta)\Bigr|\ll \frac{h}{(\log{h})^{1/10}}.
\]
\end{thrm}
\subsection{Reduce to a ternary problem}
First we use the multiplicativity of $\lambda$. Since $\lambda(n)=-\lambda(pn)$ for any prime $p$, we have
\begin{align*}
\sum_{n<x}\frac{\lambda(n)\lambda(n+1)}{n}&=\sum_{n<x}\frac{\lambda(pn)\lambda(pn+p)}{n}\\
&=\sum_{m<px}\frac{\lambda(m)\lambda(m+p)}{m}\mathbf{1}_{p|m}\\
&=\sum_{m<x}\frac{\lambda(m)\lambda(m+p)}{m}\mathbf{1}_{p|m}+O(\log{p}).
\end{align*}
Provided $p$ is small enough compared with $x$, we see that the error term is negligible (here we make important use of the logarithmic weighting). We now average this expression over $p\in [P,2P]$. This gives
\[
\sum_{n<x}\frac{\lambda(n)\lambda(n+1)}{n}=\frac{\log{P}}{P}\sum_{m<x}\frac{\lambda(m)}{m}\sum_{p\in [P,2P]}\lambda(m+p)\mathbf{1}_{p|m}+o(\log{x}).
\]
The key point of these manipulations is that we have gone from a binary problem where we want to understand the simultaneous factorizations of $n$ and $n+1$ but only have one variable at play, to a ternary problem where we want to understand the factorization of $m$ and $m+p$ and we have two variables at play. Such estimates are much more amenable to analytic estimates (this is roughly why we can solve the ternary Goldbach problem, but not the binary one). If we ignored the indicator function $\mathbf{1}_{p|m}$ and relaxed $p\in[P,2P]$ to $p\le x$ then we can estimate 
\[
\sum_{m,p<x}\lambda(m)\lambda(m+p)
\]
using the circle method. Similarly, if we ignored the indicator function $\mathbf{1}_{p|m}$ and the fact $p$ is prime, we can evaluate 
\[
\sum_{m}\Bigl|\sum_{n\in [P,2P]}\lambda(m+n)\Bigr|
\]
using the Matom\"aki-Radziwi\l\l\,Theorem. Thus our main difficulty is to replace the indicator function with something simpler, and then to control the fact we have primes in short intervals.

Since $1/m\approx 1/(m+j)$ if $|j|$ is smaller than $m$, we see that
\begin{align*}
\sum_{m<x}&\frac{\lambda(m)}{m}\sum_{p\in [P,2P]}\lambda(m+p)\mathbf{1}_{p|m}\\
&\approx\frac{1}{J}\sum_{m<x}\frac{1}{m}\sum_{j\le J}\sum_{p\in [P,2P]}\lambda(m+j)\lambda(m+j+p)\mathbf{1}_{p|m+j}.
\end{align*}
Thus the proof boils down to showing that for `most' $m\le x$ we have
\[
\sum_{j\le J}\sum_{p\in [P,2P]}\lambda(m+j)\lambda(m+j+p)\mathbf{1}_{p|m+j}=o\Bigl(\frac{J}{\log{P}}\Bigr),
\]
for a suitable choice of $J$ and $P$ of similar size, both very small compared with $x$.

\subsection{Correlation between sign patterns and residue classes}

We see that the values of $\lambda(m+j)\lambda(m+j+p)$ appearing above are determined by the vector $\mathbf{x}_m=(\lambda(m+1),\dots,\lambda(m+2P+J) )\in \{-1,1\}^{2P+J}$. Similarly, the values of $\mathbf{1}_{p|m+j}$ are determined by $m\Mod{Q}$, where $Q=\prod_{p\in[P,2P]}p$. Thus we can rewrite our expression as
\[
\sum_{\mathbf{a}\in\{-1,1\}^{2P+J}}\sum_{h\Mod{Q}}\sum_{j\le J}\sum_{p\in[P,2P]}a_ja_{j+p}\mathbf{1}_{p|h+j}\sum_{\substack{m\le x \\ m\equiv h\Mod{Q} \\ \mathbf{x}_m=\mathbf{a}}}\frac{1}{m}.
\]
In particular, if for any sign pattern $\mathbf{a}\in\{-1,1\}^{2P+J}$ we knew that the values of $m$ such that $\mathbf{x}_m=\mathbf{a}$ were equidistributed $\Mod{Q}$, then we could replace $\mathbf{1}_{p|m+j}$ with $1/p$.

Unfortunately we do not know that the values are equidistributed. However, if we temporarily fix a sign pattern $\mathbf{a}$, we see that we are left to consider
\[
\Biggl(\sum_{\substack{m\le x \\ m\equiv h\Mod{Q} \\ \mathbf{x}_m=\mathbf{a}}}\frac{1}{m}\Biggr)\cdot\sum_{p\in[P,2P]}\Bigl(\sum_{j\le J}a_ja_{j+p}\mathbf{1}_{p|h+j}\Bigr).
\]
The second expression in parentheses depends only on $h\Mod{p}$, which is independently uniformly distributed as $h$ varies $\Mod{Q}$. Therefore the sum over $p$ is a sum of many independent expressions, so we can use concentration of measure to show that for each $\mathbf{a}$ there can only be a few residue classes which cause problems.

\begin{lmm}[Application of Chernoff Bound]\label{lmm:Chernoff}
Let $\mathbf{a}\in\{-1,1\}^{2P+J}$, and let $\mathcal{E}_\mathbf{a}$ be the set of residue classes $y\Mod{Q}$ such that
\[
\Bigl|\sum_{p\in[P,2P]}\Bigl(\sum_{j\le J}a_j a_{j+p}\mathbf{1}_{p|y+j}\Bigr)-\sum_{p\in[P,2P]}\frac{1}{p}\Bigl(\sum_{j\le J}a_j a_{j+p}\Bigr)\Bigr|\ge \frac{\epsilon J}{\log{P}}.
\]
Then $\mathcal{E}_\mathbf{a}$ has size
\[
\#\mathcal{E}_{\mathbf{a}}\ll \exp\Bigl(-\frac{10P^3}{J^2\log{P}}\Bigr)Q.
\]
\end{lmm}
This means that the situation we have to rule out is that the set of $m\le x$ with $\mathbf{x}_m=\mathbf{a}$ are highly concentrated in a very small number of residue classes $\Mod{Q}$, and that for each residue class $b\Mod{Q}$ the set of $m\le x$ with $m\equiv b\Mod{Q}$ are concentrated on a very small number of sign patterns. 

\subsection{The entropy decrement argument}

Probably the most novel part of Tao's argument is handling the possibility that the residue classes are concentrated on a small set. This is the `entropy decrement argument'.

We wish to handle this possibility by showing there exists a scale (i.e. a choice of $J$ and $P$) such that the set of $m$ with $\mathbf{x}_m=\mathbf{a}$ is not concentrated in a small number of residue classes $\Mod{Q}$, for almost all relevant choices of $\mathbf{a}$. By Lemma \ref{lmm:Chernoff}, we see that this means we can effectively replace $\mathbf{1}_{p|h+j}$ with $1/p$, since the bad residue classes from $\mathcal{E}_\mathbf{a}$ have a negligible effect. 

This will be done by showing that if at some scale $(J_1,P_1)$ this is not the case, then there is a larger scale $(J_2,P_2)$ where $\mathbf{x}_m$ is concentrated in a relatively smaller set. We can repeat this observation at most a bounded number of times, and so there must be some scale $(J',P')$ when the residue classes were essentially equidistributed.

Let us first choose $(J_1,P_1)$ as our initial scales, and let 
\[
\mathbf{x}^{(1)}_m=(\lambda(m+1),\dots,\lambda(m+J_1))
\]
and $Q_1=\prod_{P_1\le p\le 2P_1}p$. Imagine for a simplification that for each residue class $b\Mod{Q_1}$, all $m\equiv b\Mod{Q_1}$ have $\mathbf{x}^{(1)}_m\in\mathcal{E}_b$ for some set $\mathcal{E}_b$ of size $\le 2^{\alpha_1 J_1}$. Let 
\[
\mathbf{x}_m^{(2)}=\Bigl(\lambda(m+1),\dots,\lambda(m+r J_1)\Bigr)
\]
be the concatenation of $\mathbf{x}_m,\mathbf{x}_{m+J_1},\dots ,\mathbf{x}_{m+(r-1)J_1}$ to form a vector at scale $J_2=r J_1$. There are $\le Q_1$ choices of a residue class for $m$ and given such a choice there are $\le 2^{\alpha_1J_1}$ choices for each $\mathbf{x}_{j}$, and so $\le Q_1 2^{\alpha_1J_2}$ choices of $\mathbf{x}_m^{(2)}$. For $r$ large, this is $2^{(\alpha_1+o(1))J_2}$.

We now let $P_2=rP_1$, and consider the relevance of residue classes modulo ${Q_2}=\prod_{P\le p\le 2P}P_2$. First imagine that for most residue classes $b\Mod{Q_2}$ the set of possible $\mathbf{x}^{(2)}_m$ with $m\equiv b\Mod{Q_2}$ is close to the maximal size $2^{\alpha_1 J_2}$ - say at least $2^{(\alpha_1-\epsilon)J_2}$. Then by double counting we see that typically at least $2^{-\epsilon J_2}Q_2$ different residue classes appear in the set of $m$ such that $\mathbf{x}_m=\mathbf{a}$ for all of the relevant choices of $\mathbf{a}$. In particular, for $\epsilon$ small enough we see this set of residue classes isn't highly concentrated. This means that at this scale we avoid the obstruction, and so we're good.

Alternatively, it could be the case that for each residue class there are many fewer possible sign patterns - say $\le 2^{\alpha_2 J_2}$ for some $\alpha_2<\alpha_1$. Then we can consider a larger scale $J_3=r_2J_2$, and the concatenation $\mathbf{x}^{(3)}_m$ of $\mathbf{x}^{(2)}_m,\mathbf{x}^{(2)}_{m+J_2},\dots,\mathbf{x}^{(2)}_{m+(r_2-1)J_2}$. There are $\le Q_2$ choices of residue class for $m\Mod{Q_2}$, and given such a choice $\le 2^{\alpha_2J_2}$ choices for each $\mathbf{x}^{(2)}_j$. Thus there are $\le Q_2 2^{\alpha_2 J_3}=2^{(\alpha_2+o(1))J_3}$ choices.

Continuing in this way (be aware that we are oversimplifying details), we either find a good scale, or a decreasing sequence $\alpha_1> \alpha_2> \dots$ such that there are $\le 2^{(\alpha_\ell+o(1)) J_\ell}$ choices of $\mathbf{x}^{(\ell)}$. By carefully controlling the quantitative aspects, we find that the $\alpha_j$ must decrease enough at each stage such that this process terminates, meaning that we are in one of the good situations above for some choice of scale $(J',P')$, although we only have very weak control over the size of $J'$ and $P'$.

The outcome of this is that we have found $J'$ and $P'$ such that for most $m\le x$
\[
\sum_{j\le J'}\sum_{p\in [P',2P']}\lambda(m+j)\lambda(m+j+p)\mathbf{1}_{p|m+j}\approx\sum_{j\le J'}\sum_{p\in [P',2P']}\frac{\lambda(m+j)\lambda(m+j+p)}{p}.
\]
\subsection{The circle method}
It suffices to show for most $m\le x$ that
\[
\Bigl|\sum_{j\le J'}\sum_{p\in[P',2P']}\lambda(m+j)\lambda(m+j+p)\Bigr|=o\Bigl(\frac{J' P'}{\log{P'}}\Bigr).
\]
As mentioned before, this can be thought of a version of the Ternary Goldbach problem in short intervals. Let $H=J'+2P'$. We see that by Fourier expansion
\begin{align*}
\sum_{j\le J'}\sum_{p\in[P',2P']}&\lambda(m+j)\lambda(m+j+p)\\
&=\frac{1}{H}\sum_{b\Mod{H}}\Bigl|\sum_{j\le H}\lambda(m+j)e\Bigl(\frac{bj}{H}\Bigr)\Bigr|^2\sum_{p\in [P',2P']}e\Bigl(\frac{-bp}{H}\Bigr),
\end{align*}
where $e(x)=e^{2\pi i x}$. The set of $b$ for which
\[
\Bigl|\sum_{p\in [P',2P']}e\Bigl(\frac{-bp}{H}\Bigr)\Bigr|=o\Bigl(\frac{P'}{\log{P'}}\Bigr)
\]
contribute a total of $o(H P'/\log{P'})$. However, the set $\mathcal{B}$ of $b$ for which the right hand side is $\gg P'/\log{P'}$ has size $O(1)$, by considering a $4^{th}$ power mean value. Thus these bad frequencies contribute
\[
\ll \frac{P'}{H\log{P'}}\sup_{b\in\mathcal{B}} \Bigl|\sum_{j\le H}\lambda(m+j)e\Bigl(\frac{bj}{H}\Bigr)\Bigr|^2.
\]
But this is $o(HP'/\log{P'})$ by the Theorem \ref{thrm:Extended}.

\section{Liouville in short intervals}\label{sec:MatomakiRadziwill}

Rather than sketch the proof of Theorem \ref{thrm:Extended} itself, we will give a sketch of just the hardest case $\theta=0$. This constitutes the main new ideas in the work of Matom\"aki and Radziwi\l\l, and is essentially just a qualitative version of Theorem \ref{thrm:MatomakiRadziwill}.
\begin{thrm}[Liouville function in almost all short intervals]
Let $\epsilon>0$. Then there is a constant $c(\epsilon)$ such that the following holds.

If $\delta X\ge c(\epsilon)$ and $0<\delta<1$, then for all but at most $\epsilon X$ values of $y\in [X,2X]$ we have
\[
\Bigl|\frac{1}{\delta y}\sum_{y\le n\le y(1+\delta)}\lambda(n)\Bigr|\le \epsilon.
\]
\end{thrm}
Informally, this shows that `almost all' intervals $[y,y(1+\delta)]$ have $\sum_{n\in[y,y(1+\delta)]}\lambda(n)=o(\delta X)$, even when $\delta X$ is going to infinity arbitrarily slowly with $X$. 
\subsection{Reduction to Dirichlet Polynomials}
To show such a result, it suffices to show
\[
\int_X^{2X}\Bigl|\frac{1}{\delta y}\sum_{y\le n\le y(1+\delta)}\lambda(n)\Bigr|^2=o(X)
\]
as $\delta X\rightarrow\infty$. Indeed, let $\mathcal{S}_\epsilon$ be the set of $y\in [X,2X]$ for which 
\[
\Bigl|\frac{1}{\delta y}\sum_{n\in[y,y(1+\delta)]}\lambda(n)\Bigr|\ge \epsilon.
\]
 Then we see that $\mathcal{S}_\epsilon$ contributes at least $\epsilon^2 m(\mathcal{S}_\epsilon)$ to the integral, where $m(\mathcal{S}_\epsilon)$ is the measure of $\mathcal{S}_\epsilon$. Therefore $m(\mathcal{S}_\epsilon)=o(X)$ for any fixed $\epsilon>0$, as required.

The standard way to attempt to estimate such integrals is to relate them to large values of Dirichlet polynomials. This can be achieved via the following lemma.

\begin{lmm}[Counting using Mellin transforms]
Let $|a_n|\le 1$ be a real sequence. Then for any $C>0$ we have
\[
\int_X^{2X}\Bigl|\frac{1}{\delta y}\sum_{y\le n\le y(1+\delta)}a_n\Bigr|^2dy \ll \frac{X}{C}+\frac{1}{X} \int_0^{C/\delta}\Bigl|\sum_{X\le n\le 3X}\frac{a_n}{n^{it}}\Bigr|^2dt .
\]
\end{lmm}

This reduces the task to showing that $\sum_n a_n/n^{it}$ cannot be large very often. In particular, by choosing $C$ as a large constant, we would be done if we could show
\begin{equation}
\int_0^{C/\delta}\Bigl|\sum_{n\le X}\frac{\lambda(n)}{n^{it}}\Bigr|^2dt=o(X^2).\label{eq:LambdaTarget}
\end{equation}
We expect \eqref{eq:LambdaTarget} to hold for $\delta X\rightarrow\infty$, but we only know how to show this when $\delta\ge X^{-5/12+\epsilon}$. This is very frustrating, since one easily gets a bound which is almost good enough from the well-known mean value theorem.

\begin{lmm}[Mean value Theorem for Dirichlet polynomials]\label{lmm:MeanValue}
Let $|a_n|\le 1$ be a complex sequence. Then
\[
\int_T^{2T}\Bigl| \sum_{n\le N}\frac{a_n}{n^{it}}\Bigr|^2dt \ll TN+N^2.
\]
\end{lmm}
This result is essentially sharp;  if $a_n=n^{it_0}$ for some fixed $t_0\in[T,2T]$ then the contribution near $t_0$ would contribute about $N^2$, whereas a random choice of $a_n$ would have a total contribution of about $TN$. The mean value estimate can be improved slightly if one is integrating over a smaller set.
\begin{lmm}[Montgomery-Halasz Mean Value Estimate]\label{lmm:Halasz}
Let $|a_n|\le 1$ be a complex sequence, and $\mathcal{U}\subseteq [T,2T]$ be a set of 1-separated points. Then
\[
\sum_{t\in \mathcal{U}}\Bigl|\sum_{n\le N} \frac{a_n}{n^{it}}\Bigr|^2\ll N^2+\#\mathcal{U}T^{1/2} N.
\]
\end{lmm}
From this, we see that one key roadblock is that it is necessary to show that $\lambda(n)$ cannot behave like $n^{it_0}$ for any $t_0\ll 1/\delta$, since this would genuinely cause problems. Fortunately, since $\lambda(n)$ is closely connected to the zeta function, we can show that $\sum_{n}\lambda(n)n^{-it}$ is always a bit smaller than the trivial bound, which at least removes this obstruction.

\begin{lmm}[Consequences of a zero-free region for $\zeta(s)$]\label{lmm:ZeroFree}
Let $|t|\le \exp((\log{x})^{5/4})$. Then there is a constant $c>0$ such that
\begin{align*}
\Bigl|\sum_{p<x}\frac{1}{p^{it}}\Bigr|&\ll \frac{x}{t\log{x}}+x\exp(-c(\log{x})^{1/7}),\\
\Bigl|\sum_{n<x}\frac{\lambda(n)}{n^{it}}\Bigr|&\ll x\exp(-c(\log{x})^{1/7}).
\end{align*}
\end{lmm}

If one happens to have coefficients which factor nicely, then the above $L^\infty$ bounds can be combined with Lemma \ref{lmm:MeanValue} to give a suitable result. For example, we see that for $(\log{P})^{100}\le T\le \exp((\log{P})^{5/4})$ we have
\begin{align*}
\int_T^{2T}\Bigl|\sum_{P\le p_1\le 2P}\frac{1}{p_1^{it}}\Bigr|^2\Bigl|\sum_{X/P\le p_2\le 2X/P}\frac{1}{p_2^{it}}\Bigr|^2dt&\ll \frac{P^2}{(\log{P})^{100}}\int_T^{2T}\Bigl|\sum_{X/P\le p_2\le 2X/P}\frac{1}{p_2^{it}}\Bigr|^2dt\\
&\ll \frac{X^2}{(\log{P})^{100}}\Bigl(1+\frac{T P}{X}\Bigr),
\end{align*}
which gives a good bound if $T<X/P$. Unfortunately one doesn't necessarily have a nice factorization like this in general. Building on the pioneering work of Vinogradov, the strongest previous work on primes and Liouville was on using combinatorial formulae to express such expressions with terms with nice factorizations, or with terms with very well-controlled coefficients.

\subsection{Factorization of typical integers}
The key idea which Matom\"aki and Radziwi\l\l\, exploit is that unlike for primes, the Dirichlet polynomial associated to the Liouville function naturally obtains many convenient factorizations because of the fact that most integers have `small' prime factors.

\begin{lmm}[Almost all integers have many small prime factors]\label{lmm:Factorize}
Let $|a_n|\le 1$ be a complex sequence. Let $\mathcal{I}_h=[\exp((\log{h})^{9/10}),h]$ with $h\le X^{1/2}$. Then for all but $O(y/(\log\log{h})^{1/4})$ values of $y\in [X,2X]$ we have
\[
\sum_{y\le n\le y(1+\delta)}a_n=\frac{10}{\log\log{h}}\sum_{y\le n\le y(1+\delta)}a_n\sum_{\substack{p|n\\ p\in\mathcal{I}_h}}1+O\Bigl(\frac{\delta X}{(\log\log{h})^{1/4}}\Bigr).
\]
\end{lmm}
Thus we can work with $a_n=\lambda(n)\sum_{p|n,p\in\mathcal{I}_h}1=-\sum_{mp=n}\lambda(m)\mathbf{1}_{p\in\mathcal{I}_h}$ instead of $\lambda(n)$. This might look like a more complicated function, but the key advantage is that now our Dirichlet polynomial can essentially be factored into two different Dirichlet polynomials, and we have control over the lengths of the two factors. After splitting the range of primes into intervals $[P,2P]$, with this trick we find it is sufficient to show that for any $P\in\mathcal{I}_h$ we have
\begin{equation}
\int_0^{C/\delta}\Bigl|\sum_{P\le p\le 2P}\frac{1}{p^{it}}\Bigr|^2\Bigl|\sum_{X/2P\le m\le 3X/P}\frac{\lambda(m)}{m^{it}}\Bigr|^2dt=o\Bigl(\frac{X^2}{(\log{P})^2}\Bigr).
\end{equation}
By Lemma \ref{lmm:ZeroFree}, we have that $\sum_m\lambda(m)/m^{it}\ll X/(P\log^{10}{X})$ if $|t|\le \log^3{X}$, which gives the result in this case.

We now consider $\log^3{X}\le |t|\le X$, so $|t|\ge \log{P}$. Provided $h\ge \exp((\log{X})^{9/10}))$, we see that $P\in\mathcal{I}_h$ implies $|t|\le X\le \exp((\log{h})^{10/9})\le \exp((\log{P})^{5/4})$. Thus, by Lemma \ref{lmm:ZeroFree}, we have 
\[
\Bigl|\sum_{P\le p\le 2P}\frac{1}{p^{it}}\Bigr|^2\ll \frac{P}{(\log{P})^2}.
\]
Combining this with the mean value estimate for $\sum_m \lambda(m)/m^{it}$ then shows this total contribution is $O(X/(\log{P})^4)$ provided that $\delta X\ge h$. Choosing $h=\exp((\log{X})^{9/10})$ gives the result for $\delta X\ge \exp((\log{X})^{9/10})$. This is already well beyond what was viewed as within reach before the Matom\"aki-Radziwi\l\l\, work, and demonstrates the power of the new idea of using such factorizations to improve on mean-value estimates. Their result goes much further, however.
\subsection{Iterative factorization}
We can use Lemma \ref{lmm:Factorize} several times for different choices of $h$ rather than just once. This gives that for a suitable choice of $h_1\le h_2\le \dots \le h_J$
\[
\sum_{y\le n\le y(1+\delta)}a_n=\prod_{i=1}^J\Bigl(\frac{10}{\log\log{h_j}}\Bigr)\sum_{y\le n\le y(1+\delta)}a_n\sum_{\substack{p_1,\dots,p_J|n\\ p_j\in\mathcal{I}_{h_j}}}1+o(\delta X).
\]
As before, we are left to show something like
\[
\int_{\log{X}}^{C/\delta}\Bigl| \sum_{\substack{X/P_1\le mp_2\dots p_J\le 2X/P_1\\ p_2\in\mathcal{I}_{h_2},\,\dots,\,p_J\in\mathcal{I}_{h_J}}}\frac{\lambda(m)}{m^{it}p_2^{it}\dots p_J^{it}}\Bigr|^2\Bigl|\sum_{P_1\le p_1\le 2P_1}\frac{1}{p_1^{it}}\Bigr|^2 dt=o\Bigl(\frac{X^2}{(\log{P_1})^2}\Bigr).
\]
The set $\mathcal{T}_1$ of $t\in[\log{X},C/\delta]$ for which $|\sum_{p_1}1/p_1^{it}|\le P_1^2/(\log{P_1})^2$ cause no problems provided $\delta X\ge h_1$. However, for the terms when $|\sum_{p_1}1/p_1^{it}|$ is large, now we cannot assume that $h_1$ is large enough that Lemma \ref{lmm:ZeroFree} can apply. We therefore have to study the exceptional set $\mathcal{E}_1=[\log{X},C/\delta]\setminus\mathcal{T}_1$. After splitting according to the size of $p_2$ we wish to show
\[
\int_{\mathcal{E}_1}\Bigl| \sum_{\substack{X/P_1P_2\le mp_3\dots p_J\le 2 X/P_1P_2\\ p_3\in\mathcal{I}_{h_3},\,\dots,\,p_J\in\mathcal{I}_{h_J}}}\frac{\lambda(m)}{m^{it}p_3^{it}\dots p_J^{it}}\Bigr|^2\Bigl|\sum_{P_2\le p_2\le 2P_2}\frac{1}{p^{it}}\Bigr|^2 dt=o\Bigl(\frac{X^2}{P_1^2(\log{P_2})^2}\Bigr).
\]
We first consider the contribution from $t\in\mathcal{T}_2$, those $t\in\mathcal{E}_1$ such that $|\sum_{p_2}p_2^{-it}|\le P_2^{0.9}$. We would still lose too much from applying the mean value theorem directly, but instead we re-insert the sum $|\sum_{p_1}p_1^{-it}|\ge P_1/(\log{P_1})^2$, to find for any $\ell>0$ the above integral is bounded by
\begin{align*}
 P_2^{1.8}\Bigl(\frac{P_1}{(\log{P_1})^2}\Bigr)^{-2\ell}\int_{\mathcal{E}_1}\Bigl| \sum_{\substack{X/P_1P_2\le mp_3\dots p_J\le 2 X/P_1P_2\\ p_3\in\mathcal{I}_{h_3},\,\dots,\,p_J\in\mathcal{I}_{h_J}}}\frac{\lambda(m)}{m^{it}p_3^{it}\dots p_J^{it}}\Bigr|^2\Bigl|\sum_{P_1\le p_1\le 2P_1}\frac{1}{p_1^{it}}\Bigr|^{2\ell} dt.
\end{align*}
If we choose $\ell$ such that $P_1^\ell\in[P_1P_2,P_1^2P_2]$, then  Lemma \ref{lmm:MeanValue} gives a bound
\[
\ll P_1^{-2\ell+o(\ell)}P_2^{1.8}\Bigl( \frac{X}{P_1P_2}\cdot (2\ell P_1)^{\ell}\Bigr)^2\ll \ell^\ell \frac{X}{P_1^{2+o(1)}P_2^{0.2}}.
\]
This is acceptable if $\ell$ is not too large, which happens if say $\log\log{P_2}\le (\log{P_1})/100$.

We are therefore left to consider the remaining set $\mathcal{E}_2=\mathcal{E}_1\setminus\mathcal{T}_2$, where both the $p_1$-sum and the $p_2$-sum are large. Continuing in this manner, we consider in turn the $p_3$-sum, the $p_4$-sum etc until the $p_J$-sum, continually splitting the integral according to whether the sum is small (in which case we can bound it) or large (in which case we consider the next sum). At the $j^{th}$ stage, we can handle the contribution from parts of the integral where
\[
\Bigr|\sum_{P_j\le p_j\le 2P_j}p_j^{-it}\Bigl|\le P_j^{1-\epsilon_j}
\]
for a suitable sequence $\epsilon_2\le \dots\le \epsilon_J\le 1/5$ (the choice $\epsilon_j=1/5-1/(5j)$ works), and  provided that the $h_j$ are also chosen to increase suitably ($h_j=\exp(h_{j-1}^{1/100j^2})$ works). Thus we are left to consider the part of the integral when all the sums are large. If $J$ is large enough, we find that $h_J\ge \exp((\log{X})^{9/10})$, and so by Lemma \ref{lmm:ZeroFree} we have 
\begin{equation}
\Bigl|\sum_{P_J\le p_J\le 2P_J}p_J^{-it}\Bigr|\le \frac{P_J}{(\log{X})^{100}}.\label{eq:FinalBound}
\end{equation}
On the other hand, we may assume $h_J\le \exp((\log{X})^{99/100})$. By considering a large power of $\sum_{p_J}p_J^{-it}$ Lemma \ref{lmm:MeanValue} shows that the set of $t\in[0,1/\delta]$ for which $|\sum_{p_J}p_J^{-it}|>P_J^{4/5}$ has cardinality $O(X^{1/2-\epsilon})$. Thus we can use \eqref{eq:FinalBound} and then apply Lemma \ref{lmm:Halasz} to deduce
\begin{align*}
\int_{\mathcal{E}_{J}}\Bigl|\sum_{P_J\le p_J\le 2P_J}\frac{1}{p_J^{it}}\Bigr|^2&\Bigl|\sum_{X/P_1\dots P_J\le m\le  2X/P_1\dots P_J}\frac{\lambda(m)}{m^{it}}\Bigr|^2dt\\
&\ll \frac{P_J^2}{(\log{X})^{200}}\Bigl(\frac{X^2}{P_1^2\dots P_J^2}+X^{1-\epsilon}\frac{X}{P_1\dots P_J}\Bigr)\\
&\ll \frac{X^2}{P_1^2\dots P_{J-1}^2(\log{X})^{200}}.
\end{align*}
Finally, this has bounded the last bit of the integral, giving the result.

For more details of the Matom\"aki-Radziwi\l\l\,Theorem, we refer the reader to the survey of Soundararajan \cite{Sound2}.

\bibliographystyle{plain}
\bibliography{Takagi}
\end{document}